\documentclass[12pt]{article}
\usepackage{amssymb}
\usepackage{amsthm}
\input{xy}
\xyoption{all}
\newcommand{\bP}{\mathbb{P}}
\newcommand{\h}{{H}_{k,b}}
\newcommand{\hh}{\overline{H}_{k,b}}
\newcommand{\mg}{\mathcal{M}_g}
\newcommand{\Pic}{\mathrm{Pic}}

\newcommand{\C}{\mathbb{C}}
\newcommand{\Q}{\mathbb{Q}}

\newtheorem{Theorem}{Theorem}
\newtheorem{Proposition}{Proposition}
\newtheorem{Lemma}{Lemma}

\title{On the rational Picard group \\ of the moduli space of curves}
\author{Claudio Fontanari}
\date{}

\begin{document}
\maketitle

\begin{abstract}
We speculate about an algebro-geometric proof of Harer's theorem on the rational 
Picard group of the moduli space of smooth complex curves. In particular,
we refine the approach of Diaz and Edidin involving the Hurwitz space 
which parameterizes smooth covers of the projective line. 
\end{abstract}

\maketitle

\section{Introduction}

The rational Picard group $\Pic(\mg) \otimes \Q$ of the moduli space $\mg$ of smooth complex curves of genus $g$ 
is unidimensional, generated by the Hodge class $\lambda$. This basic result, which turns out to be a cornerstone 
in the enumerative geometry of moduli spaces, is due to Harer \cite{Harer:83}. Indeed, according to a previous 
theorem by Mumford \cite{Mum}, $\Pic(\mg)$ can be identified with $H^2(\Gamma_g)$, where $\Gamma_g$ denotes the 
mapping class group; in \cite{Harer:83} $\Pic(\mg)$ is determined by explicitly computing $H^2(\Gamma_g)$. 
\emph{However}, as pointed out by Arbarello and Cornalba in \cite{ArbCor:08}, \emph{from the point of view of an
algebraic geometer, Harer's approach has the drawback of being entirely transcendental; in addition, his proof 
is anything but simple. It would be desirable to provide a proof of his result which is more elementary, and 
algebro-geometric in nature.} In \cite{ArbCor:08}, the trascendental part of the proof is reduced to Harer's 
computation of the cohomological dimension of $\mg$ (hence eventually to Looijenga's conjecture that $\mg$ is 
covered by $g-1$ affine open subsets) by a clever inductive procedure and a subtle spectral sequence argument. 
A different kind of reduction (namely, to Harer's stability theorem) is proposed in \cite{Tr}. Here instead 
we address the same problem by revisiting the approach of Diaz and Edidin in \cite{DE}. 

Let $\h$ be the Hurwitz space of degree $k$ covers of $\bP^1$ branched over $b$ ordered points. 
It is an \'etale cover of $(\bP^1)^b \setminus \Delta$, where $\Delta$ is the union of the large 
diagonals. Our crucial improvement on \cite{DE} consists in compatifying $(\bP^1)^b \setminus \Delta$ 
not to $(\bP^1)^b$ but to the Fulton-MacPherson space $\bP^1[b]$ (see \cite{FM}). Namely, we define the 
compactification $\hh$ of $\h$ as the normalization of $\bP^1[b]$ in the function field of $\h$.
Since the boundary of $\bP^1[b]$ is a simple normal crossing divisor, we get 

\begin{Lemma}\label{quotient}
The scheme $\hh$ has finite quotient singularities.
\end{Lemma}

On the other hand, since there is a canonical projection $\hh \to \bP^1[b] \to (\bP^1)^b$, 
we are able to adapt the construction in \cite{DE}, in particular we obtain

\begin{Proposition}\label{homology}
We have $H_1(\hh, \Q) = 0$. 
\end{Proposition}

Our main contribution is the following

\begin{Theorem}\label{main} 
We have $\dim H_2(\hh, \Q) = n$, where $n$ is the maximum number of linearly independent boundary divisors 
of $\hh$, if and only if $\Pic(\h) \otimes \Q = 0$.  
\end{Theorem} 

As a consequence of Theorem \ref{main} and \cite{DE}, Theorem 3.1 (1), Harer's result on $\Pic(\mg) \otimes \Q$ 
can be deduced in a purely algebraic way from the computation of $H_2(\hh, \Q)$, which should be approached
in the spirit of \cite{DE}, proof of Theorem~5.1~(c). We hope to come back on this in the future.

We work over the complex field $\C$. 

We are grateful to Edoardo Ballico and Gabriele Mondello for their careful reading of a previous version 
of this note.

This research has been partially supported by GNSAGA of INdAM and MIUR Cofin 2008 - 
Geo\-metria delle variet\`{a} algebriche e dei loro spazi di moduli (Italy).

\section{The proofs}

\noindent \emph{Proof of Lemma \ref{quotient}.} The induced morphism $p: \hh \to \bP^1[b]$ is a finite dominant morphism 
from a normal variety to a smooth variety and by \cite{FM}, Theorem 3, the boundary $\partial \bP^1[b]$ of $\bP^1[b]$ is 
a simple normal crossing divisor such that $p$ is smooth over $\bP^1[b] \setminus \partial \bP^1[b]$ of $\bP^1[b]$. 
Therefore by \cite{K}, Theorem~2.23, $\hh$ has finite quotient singularities. 
\qed

\noindent \emph{Proof of Proposition \ref{homology}.} Consider the canonical projection $\hh \to \bP^1[b] \to (\bP^1)^b$. 
The cellular decomposition of $(\bP^1)^b$ defined in \cite{DE}, \S 4.1, determines a cellular decomposition of
$\bP^1[b]$ by the inductive construction in \cite{FM}. This cellular decomposition lifts to $\hh$ by the proof 
of \cite{DE}, Lemma 4.1, and the corresponding cell complex does compute homology by \cite{DE}, \S 4.4. 
Hence we may argue as in \cite{DE}, proof of Theorem 5.1 (b). In particular, since there are no 1-cells in the complex, 
it follows that $H_1(\hh, \Q) = 0$, as claimed.   

\qed

\noindent \emph{Proof of Theorem \ref{main}.} Since the homology of an algebraic variety is finitely generated, 
from the Universal Coefficient Theorem for Cohomology it follows that $\dim H^2(\hh, \Q) = \dim H_2(\hh, \Q)$.
On the other hand, by Lemma~\ref{quotient} Poincar\'e duality holds for $\hh$ with rational coefficients, 
in particular we have $H_{2b-2}(\hh, \Q) \cong H^2(\hh, \Q)$. Finally, since according to the proof of 
Proposition \ref{homology} $\hh$ admits a cellular decomposition, by \cite{F}, Example 19.1.11~(b), there 
is an isomorphism $A_{b-1}(\hh) \otimes \Q \cong H_{2b-2}(\hh, \Q)$. Hence
$$
\dim A_{b-1}(\hh) \otimes \Q = \dim H_{2b-2}(\hh, \Q) = \dim H^2(\hh, \Q) = \dim H_2(\hh, \Q)
$$
and $A_{b-1}(\hh) \otimes \Q$ is generated by boundary classes if and only if 
$$
\dim H_2(\hh, \Q) = n.
$$ 
Now the conclusion follows from the exact sequence
$$
A_{b-1}(\hh \setminus \h) \to A_{b-1}(\hh) \to A_{b-1}(\h) \to 0
$$
and the equality $\Pic(\h) = A_{b-1}(\h)$.

\qed

\vspace{0.5cm}

\noindent
Claudio Fontanari \newline
Dipartimento di Matematica \newline 
Universit\`a di Trento \newline 
Via Sommarive 14 \newline 
38123 Trento, Italy. \newline
E-mail address: fontanar@science.unitn.it

\end{document}